\documentclass[a4paper]{article}
\usepackage[english]{babel}
\usepackage[utf8]{inputenc}
\usepackage{amsmath}
\usepackage{float}
\usepackage{textcomp}
\usepackage{graphicx}
\usepackage[colorinlistoftodos]{todonotes}
\usepackage{amssymb}
\usepackage{multicol}
\title{Spectrum of the Laplacian on the Vicsek Set "with no loose ends"}
\date{}
\begin{document}
\maketitle
\begin{multicols}{2}
Robert S. Strichartz

Department of Mathematics

Mallot Hall

Cornell University

Ithaca, NY, 14853

str@math.cornell.edu

\columnbreak

Sophia Zhu

Department of Computer Science

Cornell University

Ithaca, NY, 14853

jz552@cornell.edu

\end{multicols}

\underline{Abstract}: We study the spectral properties of a fractal VNLE obtained from the standard Vicsek set VS by making a countable number of identifications of points so that all the line segments in VS become circles in VNLE. We show that the standard Laplacian on VNLE satisfies spectral decimation with the same cubic renormalization polynomial as for VS, and thereby give a complete description of all eigenfunctions of the Laplacian. We then study the restrictions of eigenfunctions to the large circles in VNLE and prove that these are Lipschitz functions.

\underline{Keywords}: Vicsek set, Laplacian on fractals, spectral decimation, restrictions of eigenfunctions.

\section{Introduction}

\quad The Vicsek set (Figure 1.1) is a well-know self-similar fractal given as the attractor of an iterated function system (IFS) consisting of five similarities in the place with contraction ratios $1/3$ (see Figure 1.2).
\begin{figure}[H]
\centering
\includegraphics[width=0.6\textwidth]{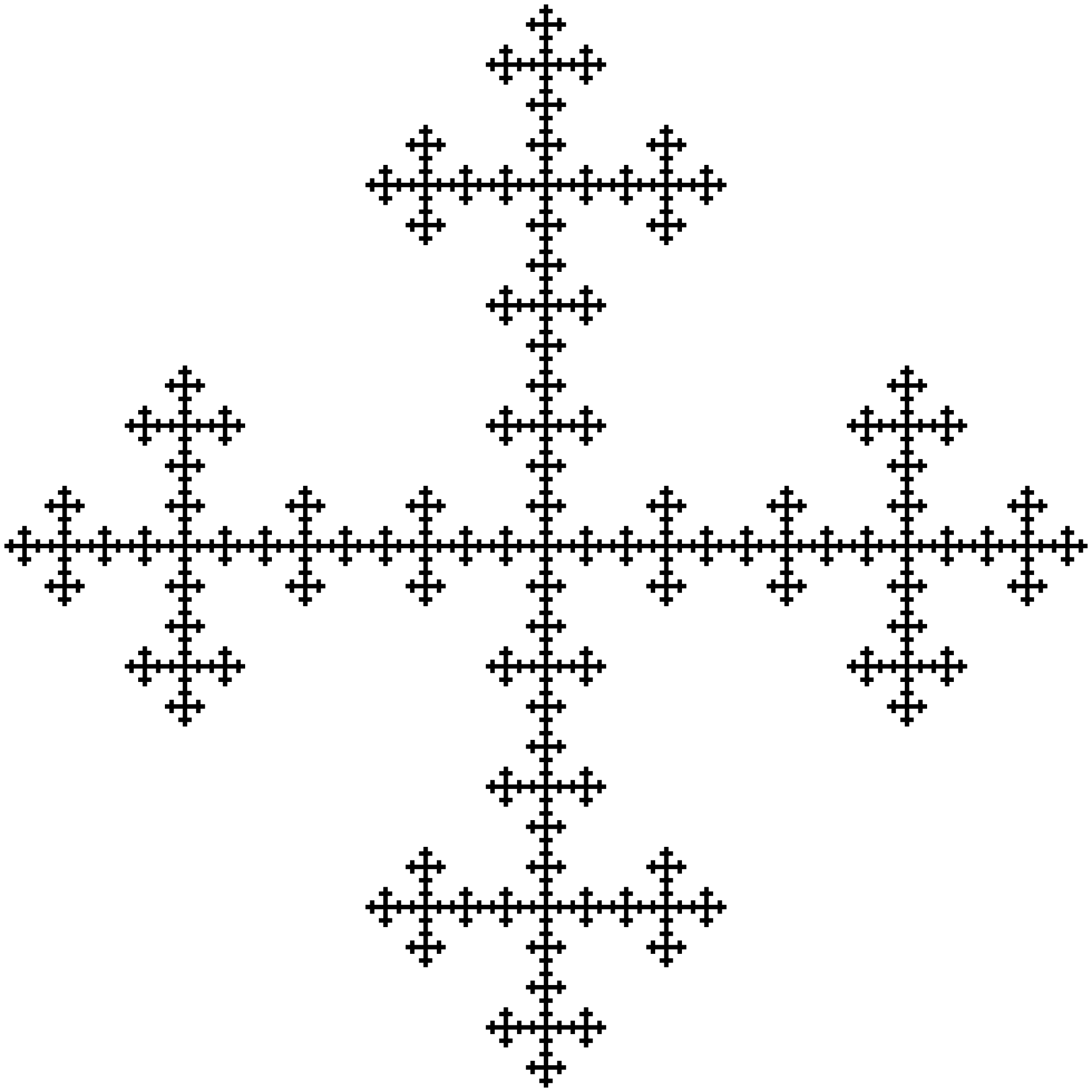}
\end{figure}
\begin{center}
Figure 1.1 : The Vicsek set.
\end{center}
\begin{figure}[H]
\centering
\includegraphics[width=0.6\textwidth]{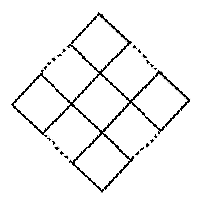}
\end{figure}
\begin{center}
Figure 1.2 : The IFS showing the five images of an initial square.
\end{center}
The Vicsek set is a example example of Kigami's class of post-critically finite (PCF) self-similar sets, and it has a completely symmetric form $\mathcal{E}$ and Laplacian $\Delta$. It was shown in $[Z]$ that the Laplacian satisfies spectral decimation, so it is possible to compute the eigenvalues and eigenfunctions exactly. See [CSU] for detail.

\quad A glance of Figure 1.1 shows that the Vicsek set has a lot of what we might call "loose ends" that lie at the end of line segments. In fact there are a countable dense set of these. These points have a local geometry that is quite different from the local geometry of generic points. It seems plausible that these points are creating a lot of "noise" in the spectrum of the Laplacian, so we would be better off without them. One way to get rid of them is to glue them together in pairs. That is what we do in this paper, producing a Vicsek set with "no loose ends", VNLE. It will turn out the VNLE is equivalent to the Julia set for a certain fifth order polynomial. The idea of introducing identifications to modify classical fractals was previously discussed in the context of turning the Sierpinski Carpet (SC) into the Magic Carpet (MC) in [BLS] and [MOS]. We believe that this will be a useful technique for lots of examples, such as the hexagasket and pentagasket, but it is not obvious how to do this explicitly.

\quad The most obvious loose ends in Figure 1.1 are the endpoints of the line segments that lie along the $x$ and $y$ axes, so in the first step of out iterative construction we identify them in pairs, the two endpoints on the x-axis will be identified, and similarly the two endpoints on y-axis. This creates two large circles that intersect at the origin. We will simultaneously define a parametrization of VNLE by a mapping g from the circle $C$ given by the unit interval with $0$ and $1$ identified. We will take $g(0)$ to be the origin, and $g$ will map $[1,\frac{1}{4}]$ to right half of the x-axis and all the points attached to it, $[\frac{1}{4},\frac{1}{2}]$ to the left half of the x-axis ending back at $g(\frac{1}{2})$ at the origin. Note that $g(\frac{1}{4})$ is the identified pair at the extremes of the x-axis. Similarly $g$ will map $[\frac{1}{2},\frac{3}{4}]$ to the top half of the y-axis and $[\frac{3}{4}, 1]$ to the bottom half, with $g(\frac{3}{4})$ the identified pair. This is the $m=0$ level approximation. Note that the parameter values $0$ and $\frac{1}{2}$ that are mapped to the same point on VNLE are mapped to the origin, not the identified pairs. At the $m=1$ level we identify the $4$ pairs of points at the ends of the next largest intervals and these correspond to the parameter values $\frac{3}{20}$, $\frac{7}{20}$, $\frac{13}{20}$, and  $\frac{17}{20}$. The parameter values $\frac{1}{5}$, $\frac{2}{5}$ both get mapped to the same center point of one of these intervals. In Figure 1.3 we show on the left the parameter circle $C$ with its identifications and VNLE with parameter values marked on the right.

\quad We iterate this procedure at all m levels, successively identifying the endpoints of the remaining longest line segments. In the limit VNLE becomes the closure of a countable union of circles. At the m-th stage of the iteration we identify the pairs of points $\frac{5^k(1 + 5j)}{2 \cdot 5^m}$ , $\frac{5^k(2 + 5j)}{2 \cdot 5^m}$ and $\frac{5^k(3 + 5j)}{2 \cdot 5^m}$ , $\frac{5^k(4 + 5j)}{2 \cdot 5^m}$ in the parameter circle as these get mapped to the intersections of intervals in $V$, while the midpoints get mapped to the loose ends in $V$ which are identified in VNLE. This description may be thought of as creating a Peano curve from $C$ to VNLE in the framework of [MOS].

\quad It is also plausible to think of this description as giving VNLE as the Julia set of a fifth degree polynomial, with the circle with identifications being the exterior ray parametrization of the Julia set. The circle is divided into five subsets. A = $[\frac{-1}{20},\frac{1}{20})$ $\cup$ $[\frac{9}{20},\frac{11}{20})$, B = $[\frac{1}{20},\frac{5}{20})$, C = $[\frac{5}{20},\frac{9}{20})$, D = $[\frac{11}{20},\frac{15}{20})$ and E = $[\frac{15}{20},\frac{19}{20})$. Each point t in the circle has a symbolic itinerary that records the sequence of subsets that ${t, 5t,5^2 t,...}$ belongs to, and points are identified exactly when they have the same itinerary. If such a polynomial exists the action of the polynomial on the Julia set is conjugate to the mapping $t \rightarrow 5t$ on the circle with identifications. Here we work entirely with the circle with identifications. See [FS, ADS, SST] for the same approach on other Julia sets.

\quad To define a Laplacian on VNLE that we denote by $K$ we consider $K$ as a limit of finite graphs $K_m$ that are obtained by taking for vertices the identified pairs on the circle at level $m$, and for edges just the intervals on the circle connecting adjacent points. Each identified pair has $4$ edges (note that some edges connect a pair to itself.) We can then define a graph Laplacian on $K_m$ by $$-\Delta_m u(x) = \sum_{z \in N_m (x)} (u(x)-u(z))$$ where $N_m (s)$ denotes the set of neighboring vertices. The Laplacian on $K$ is then the renormalized limit $$-\Delta u = \lim_{m \rightarrow \infty} (15)^m (-\Delta_m u)$$ This will be explained in detail in section 3. This definition is due to Kigami [K] (see also [S] for an elementary exposition).

\quad In this paper we give an explicit description of the spectrum of the Laplacian, both eigenvalues and eigenvectors. The spectrum of $-\Delta_m$ is described in section 2 using the method of spectral decimation of Fukushima and Shima [FS] (explained in [S]), and then the passage to the limit in section 3. Then in section 4 we describe the restrictions of eigenfunctions to circles that lie in VNLE. In particular we will prove that these are Lipschitz functions.

\quad The website [W] contains the results of numerical calculations we have perform to illustrate the results of this paper. The figures and tables in this paper show just a small sampling of these results. All eigenfunctions on VNLE are displayed as graphs on the circle with identifications.
\begin{figure}[H]
\centering
\includegraphics[width=0.9\textwidth]{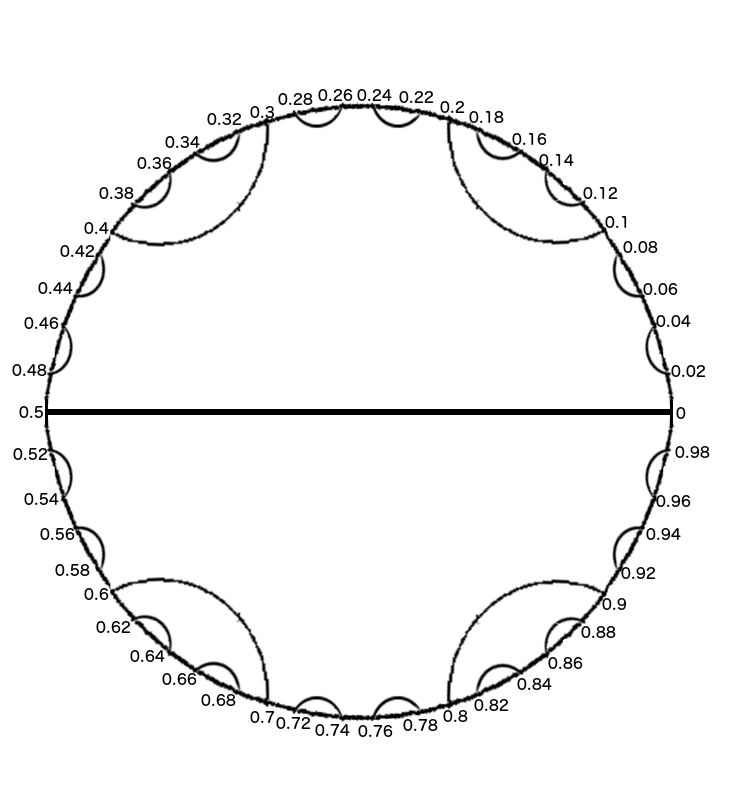}
\end{figure}
\begin{center}
Figure 1.3 : On the above graph we show the circle with identifications on level m=2. On the graph below we show an approximation to VNLE on level m=2, with the corresponding parameter values labeled.
\end{center}
\begin{figure}[H]
\centering
\includegraphics[width=0.9\textwidth]{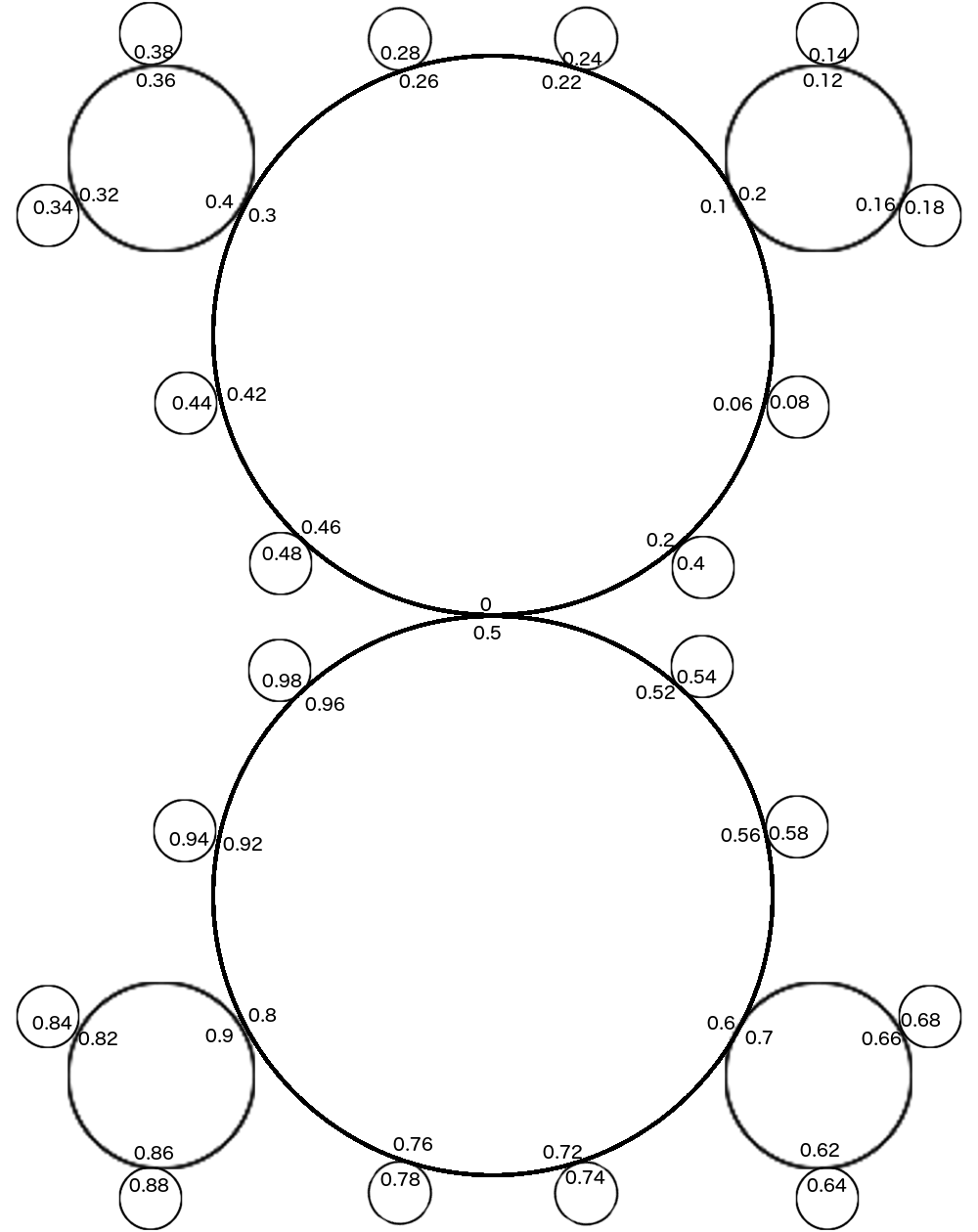}
\end{figure}

\section{Spectral Decimation}
\qquad We consider $K$ as the limit of $K_m$ as $m$ $\rightarrow$ $\infty$, where $K_m$ is a graph with $5^m$ vertices, each being an identified pair of points on the circle of the form $\frac{n}{2 \cdot 5^m}$ with $1 \leq n \leq 2 \cdot 5^m$. Note that each n can be uniquely written $n = 5^k(a + 5j)$ for a = 1,2,3,4, where $5^k$ is the highest power of $5$ dividing $n$. The identifications are given by  
\begin{flushleft}
{\large $\textnormal{(2.1)} \quad 
\frac{5^k(1 + 5j)}{2 \cdot 5^m} \sim \frac{5^k(2 + 5j)}{2 \cdot 5^m} \:\textnormal{and}\: \frac{5^k(3 + 5j)}{2 \cdot 5^m} \sim \frac{5^k(4 + 5j)}{2 \cdot 5^m}$}
\end{flushleft}
The edges of the graph $K_m$ are just the edges of consecutive points on the circle. Thus each identified pair $x \sim y$ has exactly four neighbors.
\begin{flushleft}
{\large $\textnormal{(2.2)} \quad N(x) = \{x + \frac{1}{2 \cdot 5^m}, x - \frac{1}{2 \cdot 5^m}, y + \frac{1}{2 \cdot 5^m}, y - \frac{1}{2 \cdot 5^m}\}$}
\end{flushleft}
(note that $0 = 1$ on the circle so $N(\frac{1}{2 \cdot 5^m})$ includes $1$)

\quad It is convenient to divide the vertices of $K_m$ into two types, \textit{new vertices} that are not vertices in $K_{m-1}$, and \textit{old vertices} that are vertices in $K_m$. Note that new vertices correspond to $k = 0$ and two of the neighbors will be equal to the vertex, while old vertices will have $4$ distinct neighbors not equal to the vertex. We represent $K_m$ by diagrams like Figure 2.1 for $m = 1$.
\begin{figure}[H]
\centering
\includegraphics[width=0.6\textwidth]{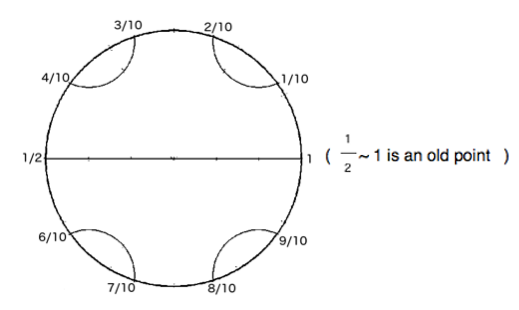}
\end{figure}
\begin{center}
Figure 2.1 : A diagram of $K_1$.
\end{center}
We represent functions on the vertices of $K_m$ by functions $u(\frac{n}{2 \cdot 5^m})$ with $u(x)=u(y)$ whenever $x \sim y$. We may then graph these functions over $[0,1]$, filling in by piecewise linear interpolation in between the points $\frac{n}{2 \cdot 5^m}$, and setting $u(0)=u(1)$.
\\
\textbf{Definition 2.1}: The Laplacian on $K_m$ is defined by 
\begin{flushleft}
{\large $\textnormal{(2.3)} \quad -\Delta_m u(x) = \underset{z \in N_m(x)}{\sum} (u(x) - u(z))$.}
\end{flushleft}
Note that for new vertices, two of the summands are $0$.
\quad We may regard $-\Delta_m$ as a symmetric $5^m \times 5^m$ matrix with diagonal entries $4$ or $2$ and off-diagonal entries $0$ or $-1$. For example, $-\Delta_1$, is 

\begin{center}
$\left (
\begin{array}{c c c c c}
2 & -1 & 0 & 0 & -1\\
-1 & 2 & 0 & 0 & -1\\
0 & 0 & 2 & -1 & -1\\
0 & 0  & -1 & 2 & -1\\
-1 & -1 & -1 & -1 & 4
\end{array}
\right )$
\end{center}

The method of spectral decimation gives a recursive description of all the eigenfunctions and eigenvalues of $-\Delta_m$. The key observation is that eigenvalues of $-\Delta_{m+1}$ either vanish identically at the old points, in which case we say that they are \textit{born} at level $m+1$, or they restrict to eigenfunctions of $-\Delta_m$ on the old points.

\quad To begin the construction we find by inspection that $-\Delta_1$ has eigenvalues $0,1,3,3,5$ with associated eigenfunctions shown in Figure 2.2. The symmetry type is indicated in parenthesis, so for example ($+$ $-$) means the function is symmetric with respect to vertical reflection and skew-symmetric with respect to horizontal reflection.
\begin{figure}[H]
\centering
\includegraphics[width=0.6\textwidth]{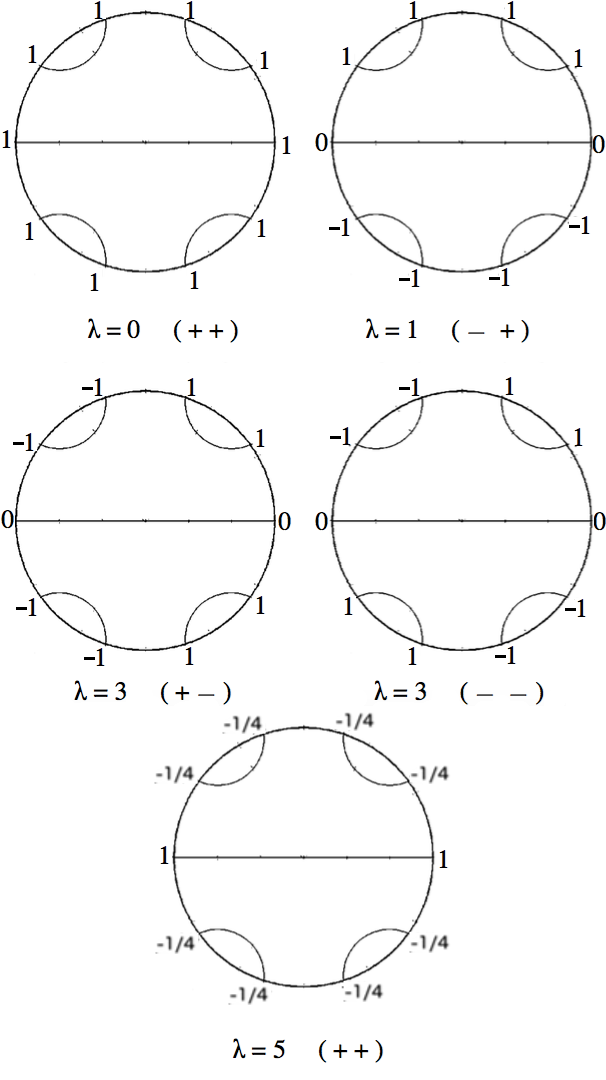}
\end{figure}
\begin{center}
Figure 2.2 : Eigenfunctions of $-\Delta_1$.
\end{center}
Since the eigenvalue $\lambda = 3$ has multiplicity $2$ we could take linear combinations of the eigenfunctions illustrated to obtain different ones supported on the upper and lower half circles. In Figure 2.3, we display the graphs of the eigenfunctions (with this choice for $\lambda = 3$).
\begin{figure}[H]
\centering
\includegraphics[width=0.8\textwidth]{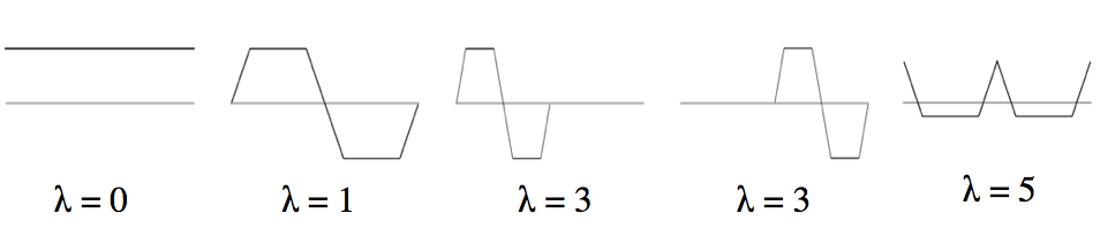}
\end{figure}
\begin{center}
Figure 2.3 : Graphs of eigenfunctions of $-\Delta_1$.
\end{center}
Note that the eigenfunctions corresponding to $\lambda = 1$ and $\lambda = 3$ are born on level $1$, while the ones corresponding to $\lambda = 0$ and $\lambda = 5$ restrict to the constant $1$ on the old vertex $\frac{1}{2} \sim 1$.

\quad In general, the eigenfunctions born on level $m+1$ will correspond to eigenvalues $\lambda = 1$ (multiplicity $5^m$) and $\lambda = 3$ (multiplicity $5^m + 1$). The other eigenfunctions will be extensions of eigenfunctions of $-\Delta_m$ on the old vertices. Generally speaking, there will be $3$ different extensions given in terms of the renormalization polynomial.
\begin{flushleft}
{\large $\textnormal{(2.4)} \quad R(x)=x(3-x)(5-x)$}
\end{flushleft}
and its three inverses $\varphi_1(x) < \varphi_2(x) < \varphi_3(x)$ on the interval $[0, 4 + \sqrt{2}]$ shown in Figure 2.4
\begin{figure}[H]
\centering
\includegraphics[width=0.6\textwidth]{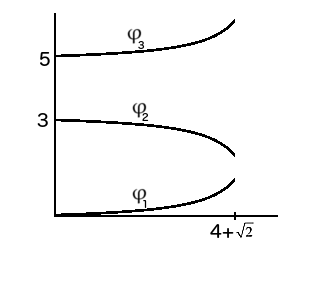}
\end{figure}
\begin{center}
Figure 2.4 : Graph of $\varphi_1$, $\varphi_2$, $\varphi_3$.
\end{center}
Note that $\varphi_1$ and $\varphi_3$ are increasing, while $\varphi_2$ is decreasing. Also $\varphi_1^{'}(0) = \frac{1}{15}$.
\\
If $\lambda$ is any eigenvalue of $-\Delta_m$ except $\lambda = 0$, then there will be $3$ extensions of the corresponding eigenfunctions from $K_m$ to $K_{m+1}$ with eigenvalues $\varphi_1(\lambda)$, $\varphi_2(\lambda)$, and $\varphi_3(\lambda)$, with the same multiplicities. Fore $\lambda=0$ there will just be $2$ extensions with eigenvalues $\varphi_1(0)=0$ and $\varphi_3(0)=5$. (Note that $\varphi_2(0)=3$ and the eigenvalue $3$ already occurs with eigenfunctions born on level $m+1$.) 
\\
There are $5^m$ eigenvalues (counting multiplicity) for $-\Delta_m$, and these extensions yield $3\cdot(5^m -1) +2=3\cdot5^m - 1$ eigenfunctions of $-\Delta_{m+1}$. Adding together the $2\cdot5^m + 1$ eigenfunctions born on level $m+1$ gives the correct number $5^{m+1}$ of eigenfunctions for $-\Lambda_{m+1}$. If 
\begin{flushleft}
{\large $\textnormal{(2.5)} \quad 0=\lambda_1^{(m)}\leq\lambda_2^{(m)}\leq...\leq\lambda_{5^m}^{(m)}$}
\end{flushleft}
denotes the enumeration of the eigenvalues of $-\Delta_m$ then 
\begin{flushleft}
$\textnormal{(2.6)} \quad 0=\varphi_1(\lambda_1^{(m)})\leq\varphi_1(\lambda_2^{(m)})\leq...\leq\varphi_1(\lambda_{5^m}^{(m)})\leq 1(multiplicity 5^m)\leq$
\end{flushleft}
\begin{flushleft}
$\varphi_2(\lambda_{5^m}^{(m)})\leq\varphi_2(\lambda_{5^m -1}^{(m)})\leq ...\leq\varphi_2(\lambda_2^{(m)})\leq 3(multiplicity \: 5^m + 1)\leq $
\end{flushleft}
\begin{flushleft}
$5 = \varphi_3(\lambda_1^{(m)})\leq\varphi_3(\lambda_2^{(m)})\leq ... \leq\varphi_3(\lambda_{5^m}^{(m)})$
\end{flushleft}

is the corresponding enumeration for $-\Delta_{m+1}$.
\\
\\
\textbf{Theorem 2.2}: Every eigenfunction of $-\Delta_{m+1}$ is either born on level $m+1$ with eigenvalue $1$ (multiplicity $5^m$) or $3$ (multiplicity $5^m + 1$), or is an extension of an eigenfunction of $-\Delta_m$ (eigenvalue $\lambda$) with eigenvalue $\lambda^{'}=\varphi_1(\lambda),\:\varphi_2(\lambda)\:\textnormal{or}\:\varphi_3(\lambda)$ (unless $\lambda = 0$, in which case only $\varphi_1(0)=0$ and $\varphi_3(0)=5$ are allowed) with the same multiplicities,
\\
\underline{Proof}: Suppose $-\Delta_m u= \lambda u$ on $K_m$. We want to find an extension $\overset{\sim}{u}$ of $u$ to $K_{m+1}$ with $-\Delta_{m+1}\overset{\sim}{u} = \lambda^{'}\overset{\sim}{u}$. We will give an algorithm for defining $\overset{\sim}{u}$ on new vertices, and then check that the eigenvalue equation holds both at new vertices and old vertices. 

\quad Note that the new vertices are sandwiched in pairs between old vertices as shown in Figure 2.5, where $A$, $B$ denote the values of $u=\overset{\sim}{u}$ at the old vertices, and $a$, $b$ denote the values of $\overset{\sim}{u}$ at the new vertices.
\begin{figure}[H]
\centering
\includegraphics[width=0.6\textwidth]{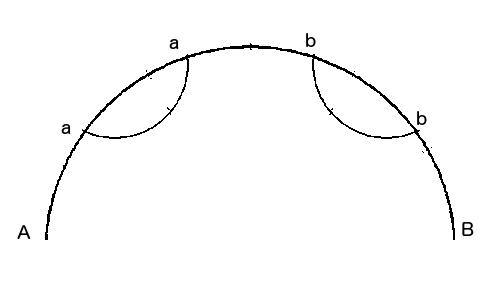}
\end{figure}
\begin{center}
Figure 2.5 : Values at old and new vertices.
\end{center}
The eigenvalue equations at the new vertices are
\begin{flushleft}
{\large $\textnormal{(2.7)} \quad 2a - A - b = \lambda^{'}a \:\textnormal{and}\: 2b - B - a = \lambda^{'}b$, }
\end{flushleft}
or equivalently 
\begin{flushleft}
{\large $\textnormal{(2.8)} \quad a + b = \frac{A + B}{1 - \lambda^{'}} \:\textnormal{and}\: a - b = \frac{A - B}{3 - \lambda^{'}}$.}
\end{flushleft}
These are easily solved to yield 
\begin{flushleft}
{\large $\textnormal{(2.9)} \quad a = \frac{(2 - \lambda^{'})A + B}{(1 - \lambda^{'})(3 - \lambda^{'})} \:\textnormal{and}\: b = \frac{A + (2- \lambda^{'})B}{(1-\lambda^{'})(3-\lambda^{'})}$.}
\end{flushleft}

Thus we take (2.9) to be out explicit extension algorithm, and this implies that the eigenvalue equation holds at new vertices.
\\
We note that $\lambda^{'} = 1$ never occurs as $\varphi_j(\lambda)$ and $\lambda^{'} = 3$ only occurs as $\varphi_2(0)$ that has been excluded, so (2.9) makes sense for all $\lambda^{'}$.

\quad We next verify that the eigenvalue equation holds at the old vertices. This part of the argument will use the fact that we began with a $\lambda$-eigenfunction of $-\Delta_m$. So consider the vertex in Figure 2.5 where the value $A$ occurs. This vertex has $4$ neighbors in $K_m$ with values denoted $B_1,B_2,B_3,B_4$ and so we have $4$ copies of Figure 2.5 with subscripts 1,2,3,4. We note that it may happen that two of these neighbors are part of the same vertex, so that $B_j = A$ (this happens at new vertices in $K_m$), but the argument is still the same. The $-\Delta_m$ eigenvalue equation is 
\begin{flushleft}
{\large $\textnormal{(2.10)} \quad \sum_{j=1}^4(A - B_j) = \lambda a$,}
\end{flushleft}
which we know is valid, and the $-\Delta_{m+1}$ eigenvalue equation would be
\begin{flushleft}
{\large $\textnormal{(2.11)} \quad \sum_{j=1}^4(A -a_j) = \lambda^{'} A$,}
\end{flushleft}
So we need to deduce (2.11) from (2.10).
\\
Note that $a_j$ are still given by (2.9), so that (2.11) is equivalent to 
$\lambda^{'}A = \sum_{j=1}^4 (A - \frac{(2-\lambda^{'})A}{(1-\lambda^{'})(3-\lambda^{'})} - \frac{B_j}{(1-\lambda^{'})(3-\lambda^{'})}) = \sum_{j=1}^4 [\frac{A-B_j}{(1-\lambda^{'})(3-\lambda^{'})} - \frac{\lambda^{'} A}{1- \lambda^{'}}]$
\\
and now using (2.10) this becomes $\lambda^{'} A = \frac{\lambda A}{(1-\lambda^{'})(3-\lambda^{'})} - \frac{4\lambda^{'}A}{1-\lambda^{'}}$. We can then cancel the factor of $A$ and simplify to obtain $\lambda = \lambda^{'}(3-\lambda^{'})(5-\lambda^{'}) = R(\lambda^{'})$. But this is the same as $\lambda^{'} = \varphi_j(\lambda)$ since $\{\varphi_j\}$ are the inverses of $R$. Thus we have verified (2.11).

\quad It remains to construct the requisite number of linearly independent eigenfunctions born on level $m+1$. We begin with eigenvalue $\lambda=3$. We consider the vertices in $K_m$ that are new in $K_j$. There are exactly $4\cdot5^{j-1}$ of these for $1\leq j\leq m$. There is also the single vertex $1 \sim \frac{1}{2}$. If $x \sim y$ is new in $K_m$ we can check that the function given by 
\begin{figure}[H]
\centering
\includegraphics[width=0.6\textwidth]{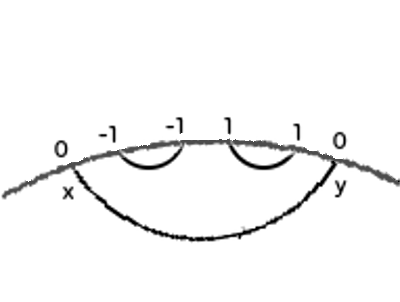}
\end{figure}
and zero elsewhere is an eigenfunction with $\lambda=3$. Similarly, if $x \sim y$ is new in $K_m$ take
\begin{figure}[H]
\centering
\includegraphics[width=0.6\textwidth]{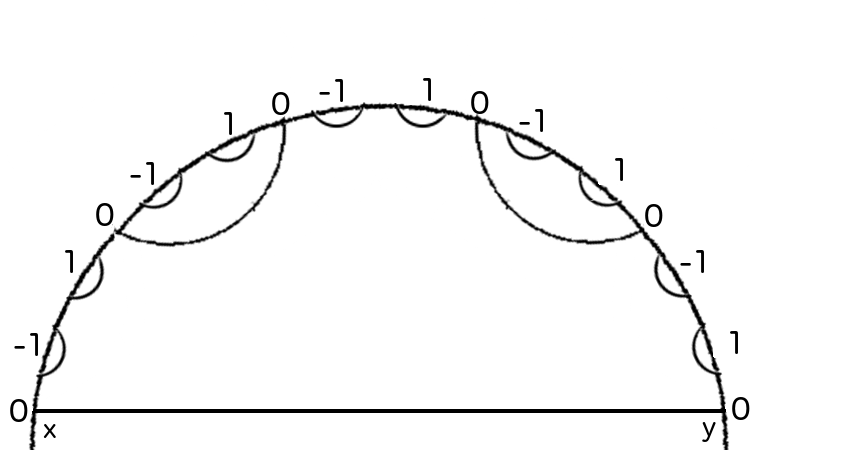}
\end{figure}
and zero elsewhere. In general, if $x \sim y$ is new in $K_j$, alternate $+1$ and $-1$ in all the new $K_{m+1}$ vertices in between $x$ and $y$. Finally, for the vertex $1 \sim \frac{1}{2}$, construct two eigenfunctions by alternating $+1$ and $-1$ at new $K_m$ vertices in either the top or bottom half of the circle. It is easy to see that these are a linearly independent set of functions and there are exactly $2 + 4(1 + 5 + 5^2 + ... + 5^{m-1}) = 5^m + 1$ of them. In terms of symmetry types, by symmetrizing or skew-symmetrizing with respect to vertical and horizontal reflections we get $(1 +5 + ... + 5^{m-1})$ of each of the $4$ symmetry types, and the $1 \sim \frac{1}{2}$ pair must have skew-symmetry with respect to the horizontal reflection giving one additional ($+$ $-$) and ($-$ $-$) type.

\quad Finally, for the eigenvalue $\lambda=1$, we note that if we consider pairs of new vertices sandwiched between two old vertices as in Figure 2.5 where now $A = B =0$, if we assign the same value to each vertex in the pair, so now $a=b$, then the 1-eigenvalue equation for $-\Delta_{m+1}$ is satisfied at the new vertices. This allows $2 \cdot 5^m$ values as variables, and the 1-eigenvalue equation at the old vertices gives $5^m$ linear constraints. It is easy to see that these constraints are linearly independent, so we are left with a space of dimension $5^m$ of 1-eigenfunctions. Although this is an indirect argument, it is easy to give examples with small support, such as 
\begin{figure}[H]
\centering
\includegraphics[width=0.6\textwidth]{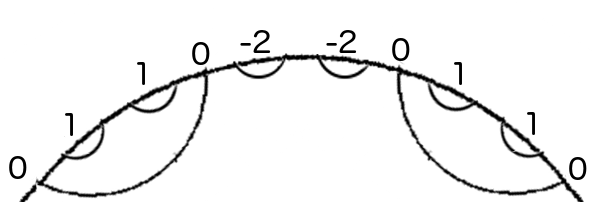}
\end{figure}
and zero elsewhere. There is also a globally supported example that alternates values $+1$ and $-1$ on consecutive pairs of new vertices, and this has symmetry type ($-$ $+$). Otherwise, all symmetry types occur with equal multiplicity.

\quad The counting argument discussed earlier shows that we have constructed $5^{m+1}$ linearly independent eigenfunctions, and hence there can be no others. $\blacksquare$

\quad \textbf{Corollary 2.3}: If we sort eigenspaces of $-\Delta_m$ into symmetry types, then the ($+$ $-$), ($-$ $+$) and ($-$ $-$) types have the same dimension, $\frac{5^m -1}{4}$, while the ($+$ $+$) type has dimension one higher.
\\
\underline{Proof}: This is obvious if we just consider the spaces of functions on $K_m$ of the different symmetry types. However it is reassuring to have an inductive argument based on our constructions. Indeed, the result is true for $m=1$ by Figure 2.2. It is clear that the extension algorithm (2.9) preserves symmetry type, so if we assume the result is true for $-\Delta_m$ then setting aside the constant function with $\lambda=0$ we obtain equal dimensions of all symmetry types in the extended eigenfunctions. The constant extends in two ways to eigenfunctions of $-\Delta_{m+1}$ with eigenvalues $0$ an $5$, both of symmetry type ($+$ $+$). However, the eigenfunctions born on level $m+1$ have one extra of each of the other three types, so again there is one extra ($+$ $+$) type. $\blacksquare$

\quad We have seen in (2.5) and (2.6) how to pass from the spectrum of $-\Delta_m$ to the spectrum of $-\Delta_{m+1}$. We insert a factor of $(15)^m$ because $\varphi_1^{'} (0) = \frac{1}{15}$, and consider the infinite sequence of spectra of $-(15)^{m}\Delta_m$ as $m \rightarrow \infty$. We claim that these converge to a limit, and in the next section we will identify the limit with the spectrum of $-\Delta$, a Laplacian on $K$.
\\
\underline{Throrem 2.4}: For any fixed $k$, we have 
\begin{flushleft}
{\large $\textnormal{(2.12)} \quad \underset{m \rightarrow \infty}{lim} (15)^m \lambda_{k} (m) = \lambda_k$}
\end{flushleft}
exists as an increasing limit (starting with $k\leq 5^m$), and 
\begin{flushleft}
{\large $\textnormal{(2.13)} \quad 0 = \lambda_0 \leq \lambda_1 \leq \lambda_2 ... $}
\end{flushleft}
\underline{Proof}: Let $m_0$ be the smallest value of $m$ satisfying $k\leq 5^m$. Then by iterating (2.6) we have 
\begin{flushleft}
{\large $\textnormal{(2.14)} \quad (15)^m \lambda_{k} (m) = (15)^{m_0} (15)^j  \varphi_1 (j) (\lambda_{k} (m_0))$}
\end{flushleft}
Where $j=m-m_0$ and $\varphi_1 (j)$ denotes the j-fold composition of $\varphi_1$. The results then follow by calculus from the property that $\varphi_1$ is $C^1$ and increasing, $\varphi_1 (0) = 0$ and $\varphi_1^{'} (0) = \frac{1}{15}$.

\begin{figure}[H]
\centering
\includegraphics[width=0.8\textwidth]{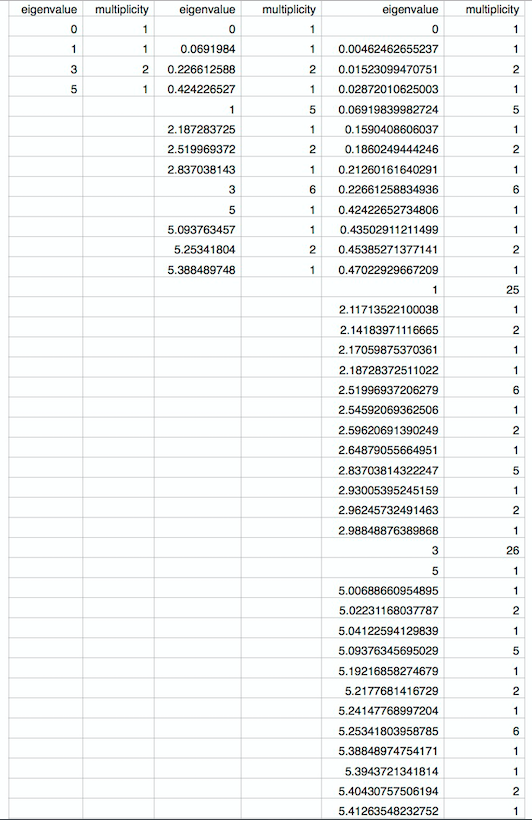}
\end{figure}
\begin{center}
Table 2.1: spectrum of $- \Delta_m$ for m = 1,2,3
\end{center}

\quad In Table 2.1 we show simultaneously the spectra of $- \Delta_m$ for m = 1,2,3. Note that all eigenvalues on level m reappear on level m+1, but not always with the same multiplicity, due to miniaturization (see section 3).

\begin{figure}[H]
\centering
\includegraphics[width=0.8\textwidth]{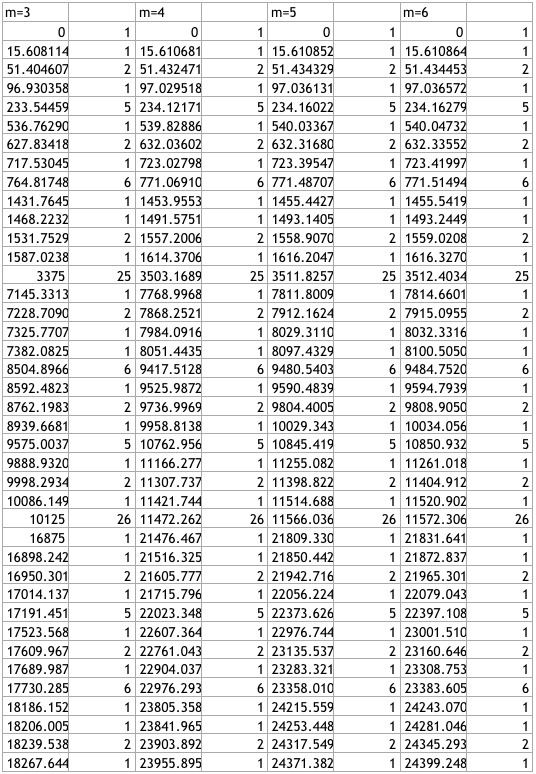}
\end{figure}
\begin{center}
Table 2.2: spectrum of $-(15)^m \Delta_m$ for m=3,4,5,6
\end{center}

\quad In Table 2.2 we show simultaneously the spectrum of $-(15)^m \Delta_m$ for m=3,4,5,6 (first 125 eigenvalue)

\quad It is clear from this table that the normalized eigenvalues converge as $m$ increases. As is the case for VS and other fractals that enjoy spectral decimation, $\lambda_{5k}(m+1) \approx 15\lambda_k(m)$, so if we define the eigenvalue counting function
\begin{flushleft}
{\large $\textnormal{(2.15)} \quad N(t) = \textnormal{number of} \{k:\lambda_k \leq t\}$}
\end{flushleft}
then $N(t)$ grows on the order of $t^\alpha$ for $\alpha = \frac{\log 15}{\log 5} = 1 + \frac{\log 3}{\log 5}$ .

\quad Thus we are led to consider the Weyl ratio
\begin{flushleft}
{\large $\textnormal{(2.16)} \quad W(t) = \frac{N(t)}{t^\alpha}$}
\end{flushleft}
As is other case for the other fractals, the Weyl ratio is asymptotically multiplicative periodic. 
\begin{flushleft}
{\large $\textnormal{(2.17)} \quad W(t) \rightarrow \varphi(t) \textnormal{as } t \rightarrow \infty \textnormal{ with}$}
\end{flushleft}
\begin{flushleft}
{\large $\textnormal{(2.18)} \quad \Phi(5t) = \varphi(t)$}
\end{flushleft}
for a discontinuous function $\phi$ that is bounded and bound away from zero. We show the graph of $N(t)$ versus $\log t$ for $-(15)^m \Delta_m$ with $m=6$ in Figure 2.6. The large jump discontinuities correspond to eigenvalues with large multiplicity, and these are responsible for the discontinuities of $\varphi$. The graph of $W(t)$ versus $\log t$ for $m=6$ is shown in Figure 2.7, and this makes the asymptotic multiplicative periodicity evident.
\begin{figure}[H]
\centering
\includegraphics[width=0.8\textwidth]{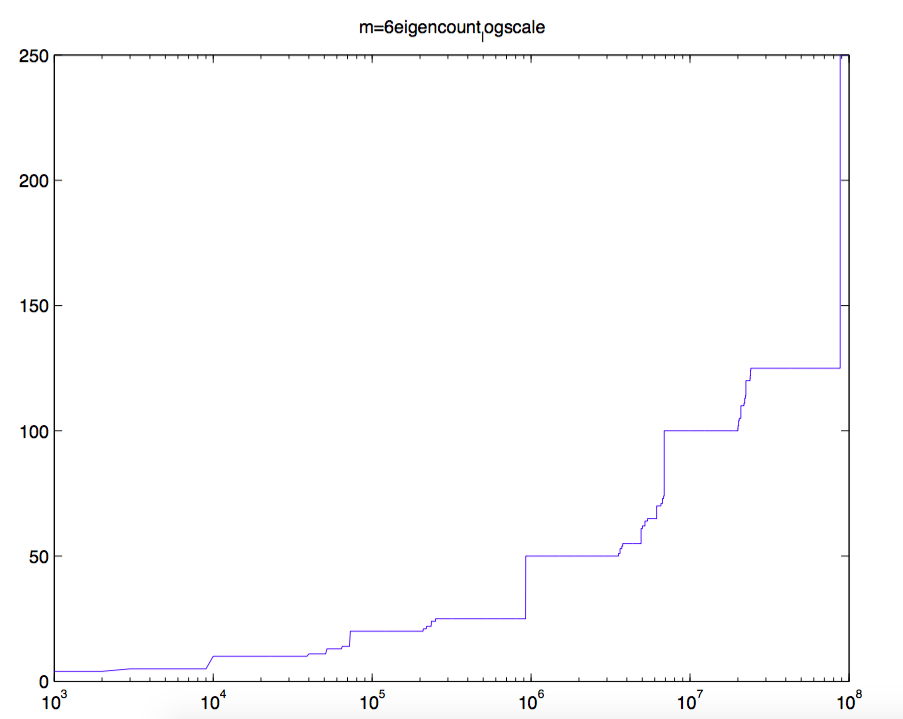}
\end{figure}
\begin{center}
Figure 2.6
\end{center}
\begin{figure}[H]
\centering
\includegraphics[width=0.8\textwidth]{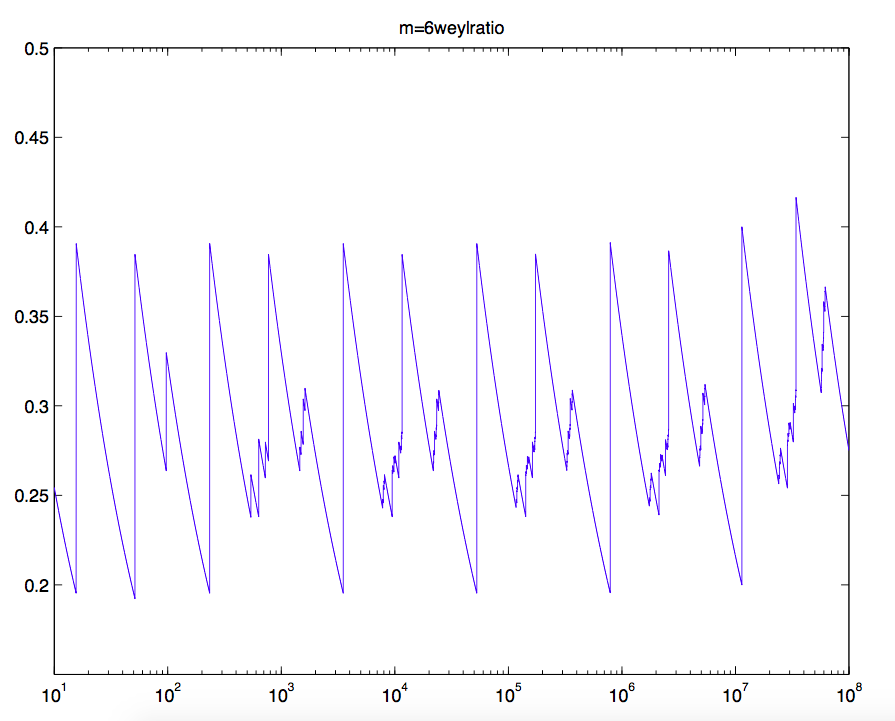}
\end{figure}
\begin{center}
Figure 2.7
\end{center}
\section{Eigenfunctions on $K$}
\qquad We now give a complete description of all eigenfunctions on $K$ as limits of eigenfunctions on $K_m$. If $u$ is a continuous function of $K$, then the restrictions of $u$ to $K_m$ determine $u$. We say that $u \in dom \Delta$ with $\Delta u = f$ for $f$ also a continuous function provided that
\begin{flushleft}
{\large $\textnormal{(3.1)} \quad \underset{m \rightarrow \infty}{lim} (15)^m \Delta_m u_{|K_m} \rightarrow f \textnormal{ uniformly on } \underset{m=1}{\overset{\infty}{\cup}} K_m$.}
\end{flushleft}
An eigenfunction of $-\Delta$ with eigenvalue $\lambda$ is a nonzero function $u \in dom \Delta$ with
\begin{flushleft}
{\large $\textnormal{(3.2)} \quad -\Delta u = \lambda u$}
\end{flushleft}
\;
\textbf{Theorem 3.1}: For every eigenfunction $u$ there exists $m_0 \geq 1$ such that $u_{|K_m}$ is a eigenfunction of $-\Delta_m$ for every $m \geq m_0$. If $\lambda_m$ denotes the eigenvalue of $u_{|K_m}$ then
\begin{flushleft}
{\large $\textnormal{(3.3)} \quad \lambda = \underset{m \rightarrow \infty}{lim} (15)^m \lambda_m$}
\end{flushleft}
where $\lambda_{m+1}=\varphi_i (\lambda_m)$ for $i = 1,2,3$ for all $m \geq m_0$. Moreover there exists $m_1$ such that $i=1$ for all $m \geq m_1$. The value of $\lambda_{m_0}$ is $0,1,3$ or $5$ if $m_0 = 1$ and $1$ or $3$ if $m_0 > 1$, and only $i=1$ or $3$ is allowed if $\lambda_m = 0$.

\underline{Proof}: Theorem 2.2 gives the construction of sequences of extensions of eigenfunctions of $-\Delta_m$ starting with $m=m_0$ with $\lambda_{m+1}=\varphi_i (\lambda_m)$. In order for the limit to exist in (3.3) we must have $\lambda_m \rightarrow 0$ as $m \rightarrow \infty$, and this implies $i=1$ for all sufficiently large $m$. Since $\varphi^{'} (0)=\frac{1}{15}$ this implies that the limit in (3.3) does exist. The extension algorithm (2.9) for $\lambda^{'} = \lambda_{m+1} = O (15^{-m})$ converges rapidly to the piecewise linear interpolation so we obtain a continuous function $u$ on $K$ in the limit, and hence an eigenfunction of $-\Delta$. $\blacksquare$

\begin{figure}[H]
\centering
\includegraphics[width=0.9\textwidth]{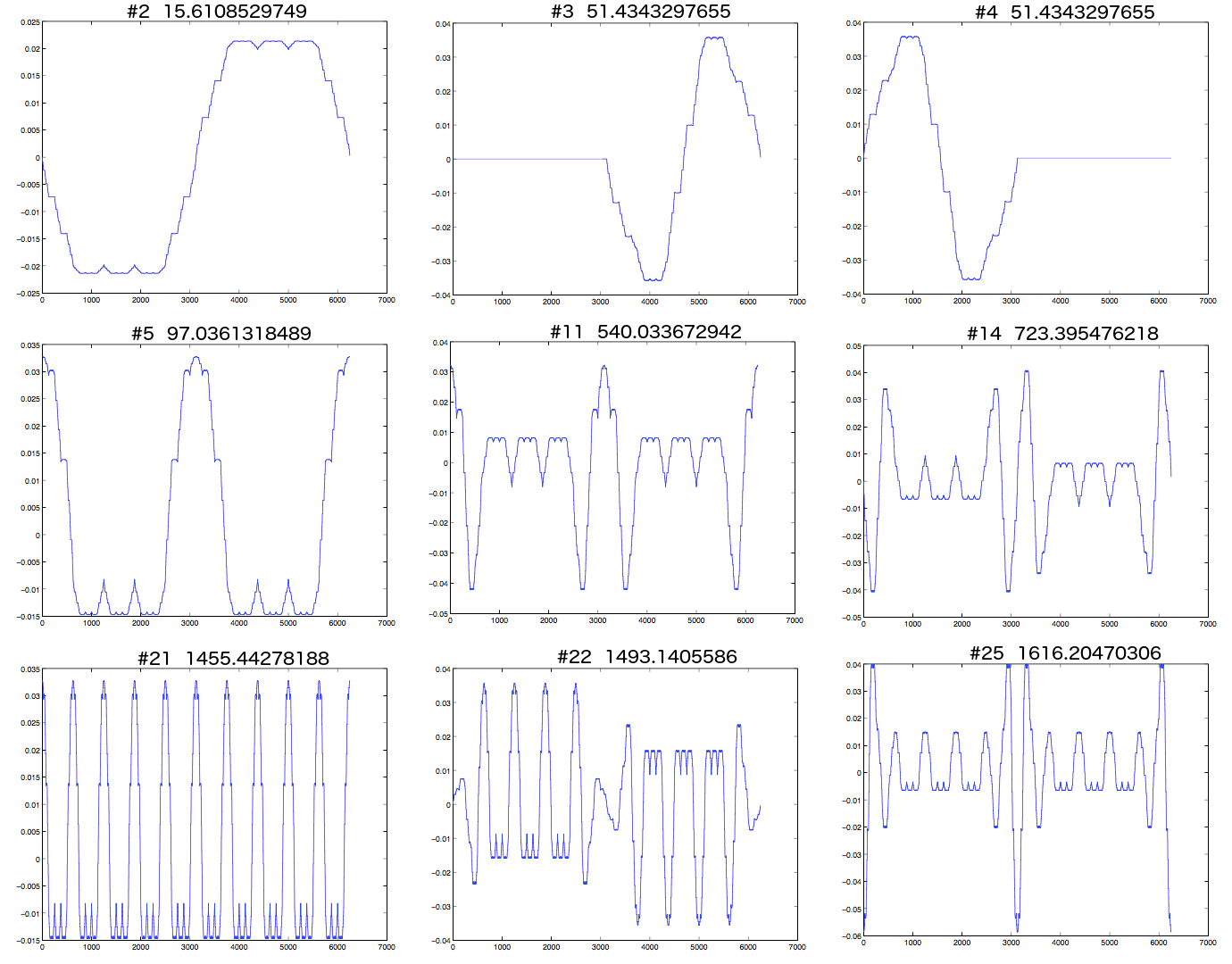}
\end{figure}
\begin{center}
Figure 3.1
\end{center}
\quad In Figure 3.1 we show a selection of the graphs of eigenfunctions. These are in fact the eigenfunctions of $-\Delta_6$, but since there is no visible change from $-\Delta_5$ to $-\Delta_6$ at the beginning of the spectrum, as can be seen on the website [w], so these are excellent approximations to the corresponding eigenfunctions of $-\Delta$. Note that 3 and 4 span a multiplicity 2 eigenspace, and 21 is a miniaturization of 5 with $u_{21} (t) = u_5 (5t)$ and $\lambda_{21} = 15 \lambda_5$. In fact it is easy to see that if $u_k$ is an eigenfunction with eigenvalue $\lambda_k$, then $u_k (5t)$ is an eigenfunction with eigenvalue $15 \lambda_k$, Thus $u_k (5t) = u_n (t)$ with $n \approx 5k$, as can be seen in Table 2.1.

\section{Restrictions of eigenfunctions}
\qquad In this section we study the restriction of eigenfunctions to the central horizontal circle $C$ in $K$. In our parametrization of $K$ this circle corresponds to a Cantor set in the upper half circle where we delete all the arcs connecting identified points. So we start with the interval $[0,\frac{1}{2}]$ (endpoints identified) and in the first stage of the construction delete $(\frac{1}{10},\frac{2}{10})$ and $(\frac{3}{10},\frac{4}{10})$ and iterate. It is convenient to define a new parameter $t$ in $[0,1]$ to describe $C$ as follows. At the first stage $0\leq t\leq \frac{1}{3}$ will parametrize $0\leq x\leq \frac{1}{10}$, $\frac{1}{3} \leq t\leq \frac{2}{3}$ will parametrize $\frac{2}{10} \leq x\leq \frac{3}{10}$, and $\frac{2}{3} \leq t\leq 1$ will parametrize $\frac{4}{10} \leq x\leq \frac{5}{10}$. And then we iterate. In terms of vertices of $K_m$ we have a simple description: the vertices of $K_m \cap C$ in consecutive order correspond to $t = \frac{k}{3^m}$ for $0\leq k \leq 3^m$ (so $k = 0$ and $k = 3^m$ are identified). This is illustrated in Figure 4.1 where the $k$ values for $K_2 \cap C$ are shown.
\begin{figure}[H]
\centering
\includegraphics[width=0.6\textwidth]{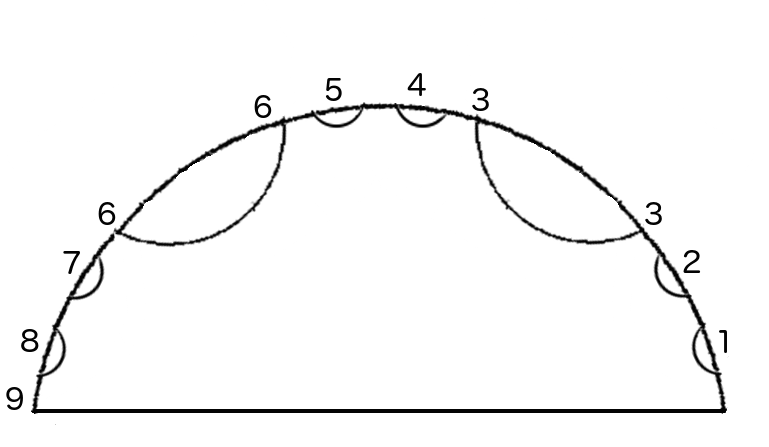}
\end{figure}
\begin{center}
Figure 4.1 : Values of $K$ for $K_2 \cap C$.
\end{center}

\quad We will only consider restrictions of primitive eigenfunctions, since all the others are miniaturizations with repetition. That means we consider eigenfunctions born on level $1$ with eigenvalue $\lambda=1,3,\textnormal{or}\:5$ and extended in all possible ways. These initial $K_1$ eigenfunctions are shown in Figure 2.2, and their restrictions to $K_1 \cap C$ as function of $t$ are shown in Figure 4.2. Note that for $\lambda=3$, which has multiplicity $2$, we only have a single restriction to $C$.
\begin{figure}[H]
\centering
\includegraphics[width=0.6\textwidth]{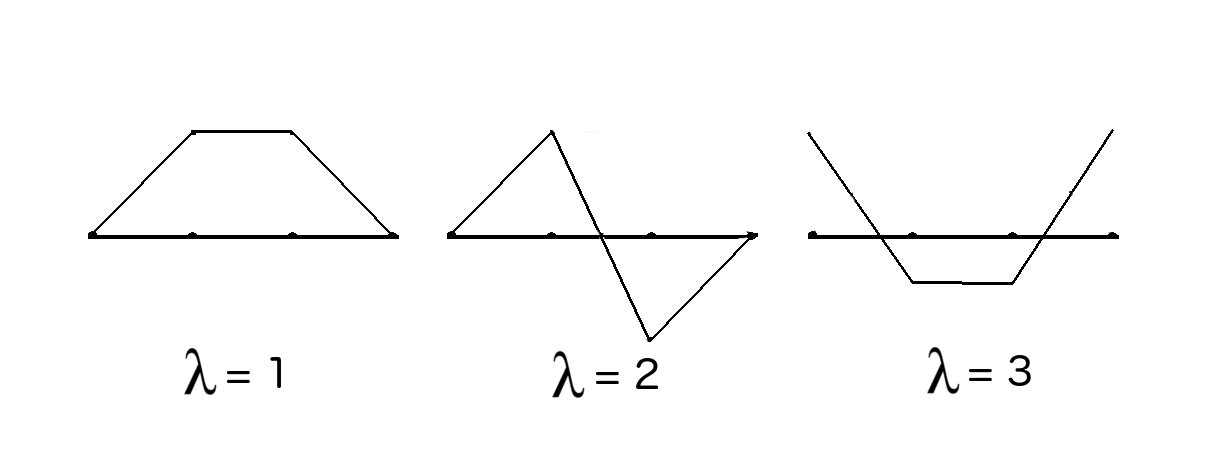}
\end{figure}
\begin{center}
Figure 4.2 : Restriction to $K_1 \cap C$.
\end{center}
We see that the $\lambda = 1$ and $\lambda=5$ restrictions are symmetric and the $\lambda=3$ restriction is skew-symmetric with respect to the reflection $t \rightarrow 1-t$. These symmetry types are preserved in the extension process.

\quad It is easy to transcribe the extension algorithm (2.9) to the restrictions. Suppose $u$ is an eigenfunction of $-\Delta_m$ with eigenvalue $\lambda$, and $f$ denotes its restriction to $K_m \cap C$, so $f(\frac{k}{3^m})$ is given. On $K_m \cap C$ the old vertices are of the form $t = \frac{3j}{3^{m+1}}$, and the new vertices are of the form $t=\frac{3j + 1}{3^{m+1}}$ and $t=\frac{3j + 2}{3^{m+1}}$, and these are conveniently sandwiched between the old vertices $\frac{3j}{3^{m+1}}$ and $\frac{3j + 3}{3^{m+1}}$. So if we choose $\lambda^{'}=\varphi_i(\lambda)$ for $i = 1,2,\textnormal{or}\:3$, then (2.9) says
\begin{flushleft}
{\large $\textnormal{(4.1)} \quad \Bigg\{\begin{array}{lcl}
f(\frac{3j + 1}{3^{m+1}}) = \frac{(2-\lambda^{'})f(\frac{3j}{3^{m+1}}) + f(\frac{3j+3}{3^{m+1}})}{(1-\lambda^{'})(3-\lambda^{'})} \\
f(\frac{3j + 2}{3^{m+1}}) = \frac{f(\frac{3j}{3^{m+1}}) + (2-\lambda^{'})f(\frac{3j+3}{3^{m+1}})}{(1-\lambda^{'})(3-\lambda^{'})}
\end{array}$}
\end{flushleft}
\;
\;
\textbf{Theorem 4.1}: let $f$ be the restriction to $C$ of an eigenfunction of $-\Delta$ on $K$. Then $f$ is Lipschitz continuous.
\\
\underline{Proof}: Let $D_m f(x) = 3^m (f(x + \frac{1}{3^m}) - f(x))$ for $x$ of the form $\frac{j}{3^m}$. It suffices to show that 
\begin{flushleft}
{\large $\textnormal{(4.2)} \quad |D_m f(x)| \leq M $}
\end{flushleft}

for some $M$ and all $m$ and $x$, for then it follows by addition that
\begin{flushleft}
{\large $\textnormal{(4.3)} \quad |f(x) - f(y)| \leq M|x - y| $}
\end{flushleft}

for $x$ and $y$ any triadic rationals, and then by continuity that (4.3) holds for all $x$ and $y$.

\quad Now we know that the eigenfunction must be obtained by extension from some level $m_0$ with all $\varphi_i$ chosen to be $\varphi_1$ for $i > m_0$. Now from (4.1) we can compute the relationship between $D_m f$ and $D_{m+1} f$, namely,
\begin{flushleft}
{\large $\textnormal{(4.4)} \quad \Bigg\{ \begin{array}{lcl}
D_{m+1} f(\frac{3j}{3^{m+1}}) = \frac{D_m f(\frac{j}{3^m})}{(1 - \lambda^{'})(1-\frac{\lambda^{'}}{3})} + \frac{3^{m+1} \lambda^{'} }{1 - \lambda^{'}}f(\frac{j}{3^m}) \\
D_{m+1} f(\frac{3j+1}{3^{m+1}}) = \frac{D_m f(\frac{j}{3^m})}{(1-\frac{\lambda^{'}}{3})} \\
D_{m+1} f(\frac{3j+2}{3^{m+1}}) = \frac{D_m f(\frac{j}{3^m})}{(1 - \lambda^{'})(1-\frac{\lambda^{'}}{3})} - \frac{3^{m+1} \lambda^{'}}{1 - \lambda^{'}} f(\frac{j+1}{3^m})
\end{array}$}
\end{flushleft}

Now since $\lambda^{'} = \varphi_1^{m + 1 -m_0} (\lambda_{m_0})$ and $\varphi_1 (\lambda) = O (\frac{1}{15} \lambda)$ we have $\frac{1}{(1-\lambda^{'})(1-\frac{\lambda^{'}}{3})} = O (\frac{1}{15^m})$ and $\frac{3^{m+1} \lambda^{'}}{1 - \lambda^{'}} = O (\frac{1}{5^m})$. So if we choose $M$ so that $|D_{m_0} f(x)| \leq \frac{M}{3}$ and $||f||_{\infty} \leq \frac{M}{3}$ a routine argument shows that (4.2) holds for $m \geq m_0$. $\blacksquare$
\\

\quad Since Lipschitz continuous functions are differentiable almost everywhere, it would be interesting to find points where we can actually compute the derivative. Because we only have access to the values $f(x)$ for $x$ a triadic rational, it seems to be unrealistic to expect to give a proof of differentiability. However, we can plausibly speculate that $f$ will have left and right derivatives at points of the form $\frac{j}{3^m}$ but not be actually differentiable there, but derivatives could exist at points of the form $\frac{j + \frac{1}{2}}{3^m}$.

\quad To be more specific, define
\begin{flushleft}
{\large $\textnormal{(4.5)} \quad a_m = \prod_{k=m}^{\infty} (\frac{1}{(1 - \frac{\lambda_k}{3})})$}
\end{flushleft}

and
\begin{flushleft}
{\large $\textnormal{(4.6)} \quad b_m = \prod_{K=m}^{\infty} (\frac{1}{(1-\lambda_k)(1 - \frac{\lambda_k}{3})})$,}
\end{flushleft}

where $\{\lambda_k\}$ is the sequence of eigenvalues of $-\Delta_k$ that correspond to the restrictions to $K_k$ of an eigenfunction. It is easy to see that the infinite products converge. By iterating the middle identity in (4.4) we obtain
\begin{flushleft}
{\large $\textnormal{(4.7)} \quad D_{m+k} f(\frac{j}{3^m} + \frac{1}{3^{m+1}} + \frac{1}{3^{m+2}} + ... + \frac{1}{3^{m+k}}) = \frac{D_m(\frac{j}{3^m})}{(\frac{1 - \lambda_{m+1}}{3})...(\frac{1 - \lambda_{m+k}}{3})}$.}
\end{flushleft}
If we knew that $f$ were differentiable at $\frac{j + \frac{1}{2}}{3^m}$ then taking the limit as $k \rightarrow \infty$ in (4.2) would yield
\begin{flushleft}
{\large $\textnormal{(4.8)} \quad f^{'} (\frac{j + \frac{1}{2}}{3^m}) = a_{m+1} D_m f(\frac{j}{3^m})$}
\end{flushleft}

\quad On the other hand, by iterating the first identity in (4.4) we obtain
\begin{flushleft}
{\large $\textnormal{(4.9)} \quad D_{m+k} f(\frac{j}{3^m}) = \frac{D_m f(\frac{j}{3^m})}{(1 - \lambda_{m+1})(1 - \frac{\lambda_{m+1}}{3}) ... (1 - \lambda_{m+k})(1 - \frac{\lambda_{m+k}}{3})} + $}
\end{flushleft}
\begin{flushleft}
{\large $\Big(\sum_{l=1}^k \frac{3^{m+l} \lambda_{m+l} (1 - \frac{\lambda_{m+l}}{3})}{(1 - \lambda_{m+l})(1 - \frac{\lambda_{m+l}}{3}) ... (1 - \lambda_{m+k})(1 - \frac{\lambda_{m+k}}{3})} \Big)f(\frac{j}{3^m})$.}
\end{flushleft}

If we knew that the right derivative of f existed at $\frac{j}{3^m}$ then we could take the limit in (4.9) to obtain the expression.
\begin{flushleft}
{\large $\textnormal{(4.10)} \quad b_{m+1} D_m f(\frac{j}{3^m}) + (\sum_{k=1}^{\infty} 3^{m+k} \lambda_{m+k} (1 - \frac{\lambda_{m+l}}{3}) b_{m+k})f(\frac{j}{3^m})$}
\end{flushleft}
for the right derivative. Similarly, using the bottom identity in (4.4) we would obtain the expression
\begin{flushleft}
{\large $\textnormal{(4.11)} \quad b_{m+1} D_m f(\frac{j-1}{3^m}) - (\sum_{k=1}^{\infty} 3^{m+k} \lambda_{m+k} (1 - \frac{\lambda_{m+l}}{3}) b_{m+k})f(\frac{j-1}{3^m})$}
\end{flushleft}
for the left derivative. It is easy to see that the infinite sums in (4.10) and (4.11) converge. There seems to be no reason why (4.10) and (4.11) should be equal, so we conjecture that the full derivative does not exist at these triadic rational points.

\quad The above arguments are based on the observation that (4.1) is rather close to piecewise linear interpolation between old vertices and new vertices when $\lambda^{'}= \varphi_1 (\lambda)$. Especially if we go to the region where $\lambda_{m+1} = \varphi_1 (\lambda_m)$ for all $m \geq m_0$, the correction to piecewise linear interpolation goes to zero very rapidly. This becomes very apparent in the graphs shown in Figure 4.3. In a similar vein we can describe what happens when $\lambda^{'}= \varphi_2 (\lambda)$ or $\lambda^{'}= \varphi_3 (\lambda)$.

\quad If $\lambda^{'}= \varphi_3 (\lambda)$ with $\lambda$ close to zero then $\lambda^{'}$ is close to $5$ and (4.1) says $f(\frac{3j + 1}{3^{m+1}}) \approx -\frac{3}{8} f (\frac{j}{3^m}) + \frac{1}{8} f (\frac{j+1}{3^m})$ and $f(\frac{3j + 2}{3^{m+1}}) \approx \frac{1}{8} f (\frac{j}{3^m}) - \frac{3}{8} f (\frac{j+1}{3^m})$. In particular, if $f(\frac{j}{3^m})$ and $f(\frac{j+1}{3^m})$ are of the same sign and approximately equal, then $f(\frac{j+1}{3^m})$ and $f(\frac{j+2}{3^m})$ are of the opposite sign and approximately equal to $-\frac{1}{4} f(\frac{j}{3^m})$. This is apparent in Figure 4.4.

\quad If $\lambda^{'}=\varphi_2 (\lambda)$ with $\lambda$ close to zero then $\lambda^{'} \approx 3- \frac{\lambda}{6}$ and (4.1) says $f(\frac{3j+1}{3^{m+1}}) \approx -f(\frac{3j+2}{3^{m+1}}) \approx \frac{3}{\lambda}(f(\frac{j}{3^{m}}) - f(\frac{j+1}{3^{m}}))$. This creates very large fluctuations at new vertices between old vertices where $f$ is moderately varying, but these disappear when $f$ is almost constant on old vertices . This behavior can be seen in Figure 4.5.

\begin{figure}[H]
\centering
\includegraphics[width=1.0\textwidth]{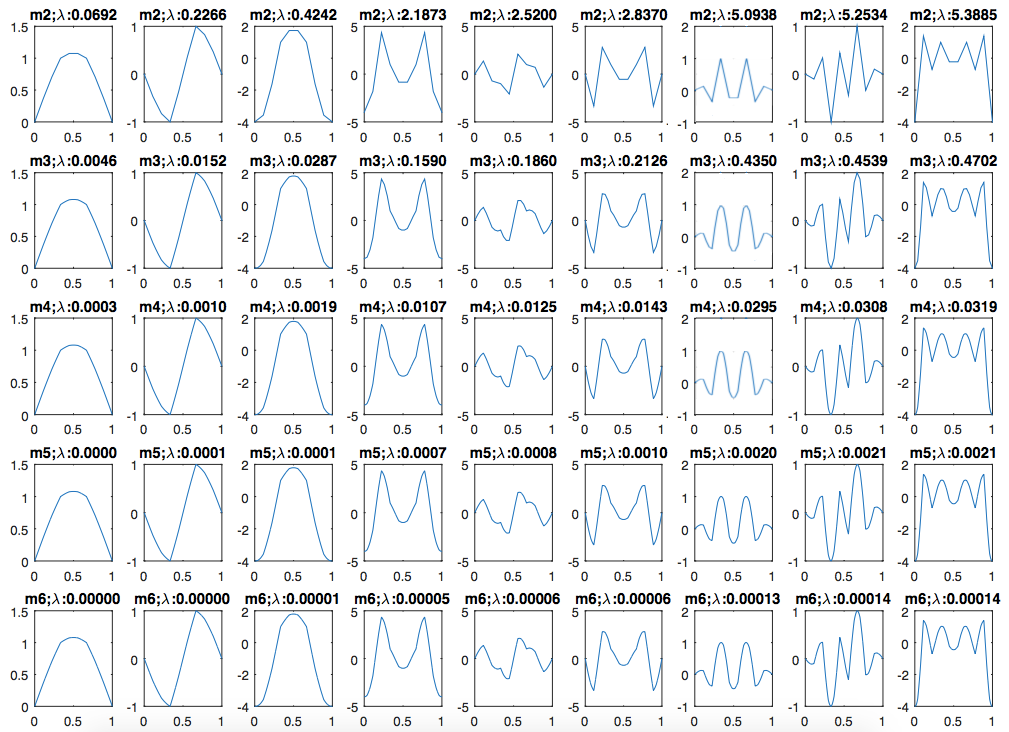}
\end{figure}
\begin{center}
Figure 4.3 : The first nine restrictions of nonconstant eigenfunctions for $m=2,3,4,5,6$ Each column corresponds to the choice $\lambda_{m+1} = \varphi_1 (\lambda_m)$.
\end{center}

\begin{figure}[H]
\centering
\includegraphics[width=0.65\textwidth]{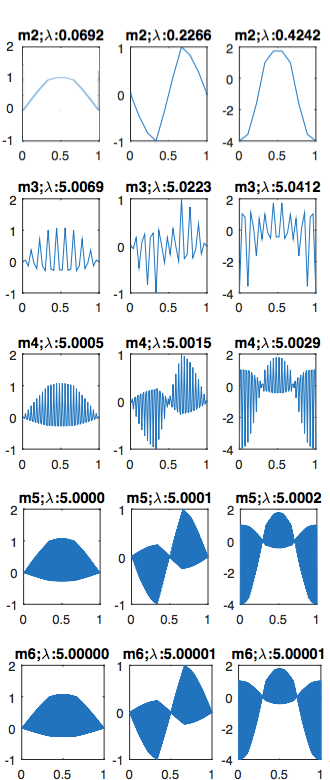}
\end{figure}
\begin{center}
Figure 4.4 : Restrictions of eigenfunctions corresponding to the choices $\lambda_m = \varphi_3(\lambda_{m-1})$ and $\lambda_k = \varphi_1(\lambda_{k-1})$ for $2 \leq k < m$, starting with the first three noncontant eigenfunctions in each column.
\end{center}

\begin{figure}[H]
\centering
\includegraphics[width=0.8\textwidth]{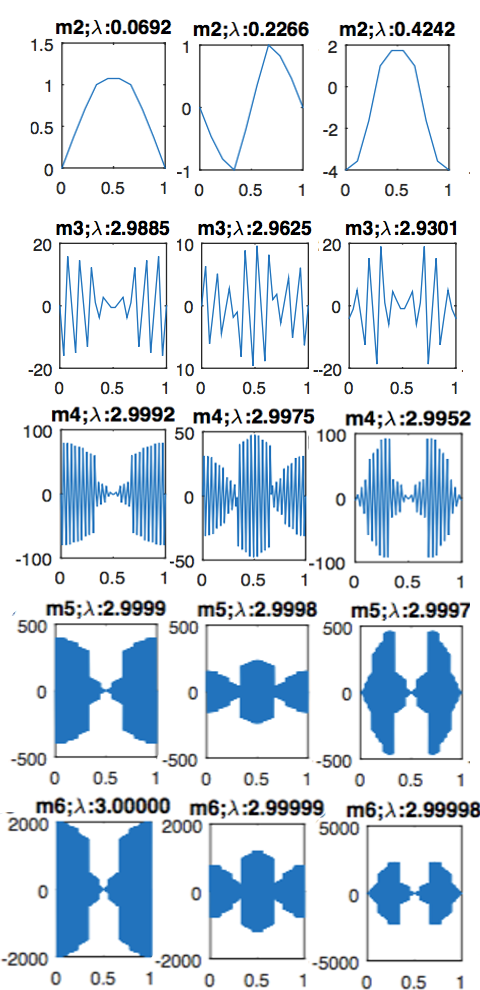}
\end{figure}
\begin{center}
Figure 4.5 : Restrictions of eigenfunctions corresponding to the choices $\lambda_m = \varphi_2(\lambda_{m-1})$ and $\lambda_k = \varphi_1(\lambda_{k-1})$ for $2 \leq k < m$, starting with the first three noncontant eigenfunctions in each column. Note the stretching of the vertical scale as m increases. 
\end{center}
\begin{center}
\textbf{References}
\end{center}
[ADS] T. Aougab, S. C. Dong and R. S. Strichartz, Laplacians on a family on quadratics Julia sets II, \textit{Commun. Pure Appl. Anal.} 12(2013) 1-58. \break
[BLS] J. Bello, Y. Li and R. S. Strichartz, Hodge-deRham theory of K-forms on carpet type fractals, in \textit{Excursions in Harmonic Analysis}, Vol.3. (eds.) R. Balan, M. J. Begue, J. J. Benedetto, W. Czaja and K. A. Okoudjou (Birkhauser, 2015), pp. 23-62. \break
[CSW] S. Constantin, R. Strichartz and M. Wheeler, \textit{Analysis of the Laplacian and spectral operators on the Vicsek set}, Comm. Pure Appl. Aoal, 10(2011), 1-44.\break
[FS] M. Fukushima and T. Shima. On a spectral analysis for the Sierpinski gasket. \textit{Potential Anal.}, 1(1):1-35, 1992. \break
[FS] T. Flock and R, Strichartz, Laplacians on a family of quadratic Julia sets, \textit{Trans. Am. Math. Soc.} 264(8) (2012) 3915-3965. \break
[K] Jun Kigami, "Analysis on Fractals," volume 143 of Cambridge Tracts in Mathematics, Cambridge University Press, Cambridge, 2001. \break
[MOS] Denali Molitor, Nadia Ott and Robert Strichartz, \textit(Using Peano Curves to Construct Laplacians on Fractals), Fractals, Vol. 23, No. 4(2015) 1550048 (29 pages) \break
[RT] Luke G. Rogers and Alexander Teplyaev, \textit{Laplacians on the basilica Julia set}, Comm. Pure Appl. Anal. 9(2010), 211-231. \break
[S] Robert S. Strichartz, "Differential Equations on Fractals, A Tutorial," Princeton University Press, Princeton, NJ, 2006. \break
[SST] C. Spicer, R. S. Strichartz and E. Totari, Laplacians on Julia Sets III: Cubic Julia sets and formal matings, \textit{Contemp. Math.} 600 (2013) 327-348. \break
[W] http://www.math.cornell.edu/$\sim$jz552 \break
[Z] D. Zhou. Spectral analysis of Laplacians on the Vicsek set. \textit{Pacific J. Math.}, 241(2):369-398, 2009.
\end{document}